\documentclass[11pt,twoside]{amsart}
\usepackage{amsmath}
\usepackage{amsthm}
\usepackage{amsfonts}
\usepackage{amssymb}
\usepackage{latexsym}
\usepackage[all]{xy}
\usepackage{mathrsfs}
\usepackage{bbm}
\usepackage{enumerate}
\usepackage{tikz-cd}
\usepackage{hyperref}
\numberwithin{equation}{section} \theoremstyle{plain}
\newtheorem*{thm*}{Main Theorem}
\newtheorem{thm}{Theorem}[section]

\newtheorem*{cor*}{Corollary}
\newtheorem{lem}[thm]{Lemma}
\newtheorem*{lem*}{Lemma}
\newtheorem{prop}[thm]{Proposition}
\newtheorem*{prop*}{Proposition}

\newtheorem*{rem*}{Remark}

\newtheorem*{exam*}{Example}
\newtheorem{defn}[thm]{Definition}
\newtheorem*{defn*}{Definition}

\newtheorem*{conj*}{Conjecture}
\newtheorem*{ack*}{ACKNOWLEDGEMENTS}

\DeclareMathOperator{\Irr}{Irr}

\newcommand{\Gkdim}{\mbox{\rm GKdim}}
\newcommand{\GL}{\mbox{\rm GL}}

\newcommand{\id}{\mbox{\rm id}}

\newcommand{\YD}{\mbox{$\mathcal{YD}$}}

\newcommand{\Bb}{\mbox{$\mathcal{B}$}}

\newcommand{\Oo}{\mbox{$\mathcal{O}$}}
\newcommand{\CC}{\mbox{$\mathbb{C}$}}

\newcommand{\GG}{\mbox{$\mathbb{G}$}}
\newcommand{\NN}{\mbox{$\mathbb{N}$}}

\newcommand{\ZZ}{\mbox{$\mathbb{Z}$}}
\newcommand{\kk}{\mbox{$\mathbbm{k}$}}

\begin{document}

\title{Finite GK-dimensional Nichols Algebras over the  quaternion group $\mathbb{Q}_8$ }
\author{Yongliang Zhang}
\address{School of Mathematics and Statistics, Ankang University,
Ankang 725000, China}
\email{yongliang\_zhang87@yeah.net}
\date{}

\begin{abstract}
 We contribute to the classification of Hopf algebras with finite
 Gelfand-Kirillov dimension, GK-dimension for short, through the study of 
 Nichols algebras over $\mathbb{Q}_8$,the  quaternion group . 
 We find all the irreducible Yetter-Drinfeld modules $V$ over $\mathbb{Q}_8$, and determine which Nichols algebras $\Bb(V)$ of $V$ are finite GK-dimensional, and 
which Nichols algebras $\Bb(V)$ of $V$ are finite dimensional.
 \end{abstract}
\maketitle

\tableofcontents
\section{\bf Introduction}
For non-commutative algebras,Gelfand-Kirillov dimension (GK-dimension for short) is an important tool as well as the Krull dimension. In recent years, Hopf algebras with finite Gelfand-Kirillov dimension  received considerable attention, see \cite{A5,A6,W1,Wu,Z0,Z}. By the lifting method appearing
in \cite{A2}, the author considered the following problems  in \cite{A1}:
\begin{enumerate}
\item When is $\dim \Bb(V)$ finite?
\item When is $\Gkdim \Bb(V)$ finite?
\end{enumerate}
For the first problem, see \cite{A3} when $\Bb(V)$ is of diagonal type.  Let $G$ be a group, and $\CC$ the field of all complex numbers. One problem is to find
all the Nichols algebras $\Bb(V)$ with finite dimension 
for any $V\in{ _G^G}\YD$, the Yetter-Drinfeld modules over the group algebra $\kk G$. The cases when $G$ is a finite
simple group were studied in \cite{A4,A8,A9,A10,A11,A12}.
For the second problem, great progress was acheived when 
$G$ is an abelian group, see \cite{A5,A6}. For non-abelian groups, the author considered the infinite dihedral group $\mathbb{D}_{\infty}$ in \cite{Z0}, and classified the Nichols algebras with finite GK-dimension.
In this paper we deal with the quaternion group $\mathbb{Q}_{8}$, and aim to find all the irreducible Yetter-Drinfeld modules $V$ over $\mathbb{Q}_{8}$. Among the irreducible modules $V$  we classify those $V$ with $\Gkdim\Bb(V)<\infty$, and those $V$ with $\dim\Bb(V)<\infty$.
Precisely, we prove the main theorems.
\begin{thm}
The only Nichols algebras of the finite dimensional irreducible Yetter-Drinfeld
 modules over $\kk \mathbb{Q}_{8}$ with finite GK-dimension,
 up to isomorphism, are those in the following list.
\begin{enumerate}[\rm (1)]
\item $\Bb(\Oo_1,\rho_i)$, with $1\leq i\leq 5$.
\item $\Bb(\Oo_{x},\phi_{0})$, $\Bb(\Oo_{x},\phi_{2})$
\item $\Bb(\Oo_{x^2},\rho_{i})$, with $1\leq i\leq 5$.
\item $\Bb(\Oo_{y},\phi_{0})$, $\Bb(\Oo_{y},\phi_{2})$.
\item $\Bb(\Oo_{xy},\phi_{0})$, $\Bb(\Oo_{xy},\phi_{2})$.
\end{enumerate}
\end{thm}

\begin{thm}
The only Nichols algebras of the finite dimensional irreducible Yetter-Drinfeld
 modules over $\kk \mathbb{Q}_{8}$ with finite dimension,
 up to isomorphism, are those in the following list.
\begin{enumerate}[\rm (1)]
\item  $\Bb(\Oo_{x},\phi_{2})$
\item $\Bb(\Oo_{x^2},\rho_{5})$.
\item  $\Bb(\Oo_{y},\phi_{2})$.
\item $\Bb(\Oo_{xy},\phi_{2})$.
\end{enumerate}
\end{thm}
For the proof of these theorems, see Proposition \ref{prop:3.1},\ref{prop:4.1},\ref{prop:5.1},\ref{prop:6.1} and Proposition\ref{prop:7.1}.
\section{\bf Preliminaries}
\subsection{Notations}
Let $\NN$ and $\ZZ$  be the set of positive integers and the set of all integers, respectively. Let $G$ be a group and $\kk$ an algebraic closed field with characteristic zero. In particular, we take $\kk=\CC$, the field
of all complex numbers. Let $\mathbf{i}=\sqrt{-1}$ be the imaginary unit, and $q=e^{\pi \mathbf{i}/2}\in \GG_4$. For any $g,h\in G$, we write $g\rhd h=ghg^{-1}$, for the adjoint action of $G$ on itself.
The conjugacy class of $g$ in $G$ will be denoted by $\Oo_g$,
and the centralizer of $g$ is denoted by $G^g=\{ h\in G|hg=gh\}$. 
\subsection{Yetter-Drinfeld modules and Nichols algebras}
Let $V$ be a vector space and $c\in \GL(V\otimes V)$. Then the pair
$(V, c)$ is said to be a braided vector space if $c$ is a solution of 
the braid equation
\begin{equation*}
(c\otimes \id)(\id\otimes c)(c\otimes \id)=(\id\otimes c)(c\otimes \id)(\id\otimes c).
\end{equation*}
We call a braided space $(V, c)$ diagonal type if there exists a matrix $\mathfrak{q}=(q_{ij})_{i,j\in \mathbb{I}}$ with $q_{ij}\in \kk^{\times}$ and $q_{ii}\neq 1$ for any $i,j\in \mathbb{I}$ such that
\begin{equation*}
c(x_i\otimes x_j)=q_{ij}x_j\otimes x_i, i,j\in \mathbb{I}.
\end{equation*}
\begin{defn}\cite[Definition 3]{A1}
Let $H$ be a Hopf algebra with bijective antipode $S$. 
A Yetter-Drinfeld
module over $H$ is a vector space $V$ provided with
\begin{enumerate}[$(1)$]
\item a structure of left $H$-module $\mu: H\otimes V\to V$ and
\item a structure of left $H$-comodule $\delta: V\to H\otimes V$, such that
\end{enumerate}
for all $h\in H$ and $v\in V$, the following compatibility condition holds:
\begin{align*}
\delta(h\cdot v)=h_{(1)}v_{(-1)}S(h_{(3)})\otimes h_{(2)}\cdot v_{(0)}.
\end{align*}
The category of left Yetter-Drinfeld modules is denoted by $_H^H\YD$.
\end{defn}
In particular, if $H=\kk G$ is the group algebra of the group $G$, then a Yetter-Drinfeld module over $H$ is a $G$-graded vector space $M=\oplus_{g\in G}M_g$ provided with a $G$-module structure such that $g\cdot M_h=M_{ghg^{-1}}$.

It can be shown that each Yetter-Drinfeld module $V\in {_H^H\YD}$
is a braided vector space with braiding
\begin{align*}
c(x\otimes y)=x_{(-1)}\cdot y\otimes x_{(0)},\quad x,y\in {_H^H\YD}.
\end{align*}

The category of Yetter-Drinfeld modules over $\kk G$ is denoted by $_G^G\YD$. Let $\Oo\subseteq G$ be a conjugacy class of $G$, then we denote by ${_G^G}\YD(\Oo)$ 
the  subcategory of ${_G^G}\YD$ consisting of all $M\in {_G^G}\YD$ with $M=\bigoplus_{s\in \mathcal{O}}M_s$.
\begin{defn}\cite[Definition 1.4.15]{H1}
Let $g\in G$, and let $V$ be a left $\kk G$-module. Define
\begin{align*}
M(g,V)=\kk G\otimes_{\kk G^g}V
\end{align*}
as an object in $_G^G\YD(\Oo_g)$, where $M(g,V)$ is the induced $\kk G$-module, the $G$-grading is given by 
\begin{align*}
\deg (h\otimes v)=h\rhd g, \quad \text{for all}~h\in G, v\in V,
\end{align*}
and the $\kk G$-comodule structure is 
\begin{align*}
\delta(h\otimes v)=(h\rhd g)\otimes (h\otimes v).
\end{align*}
\end{defn}
Let $V\in {_G^G}\YD$. Let $I(V)$ be the largest coideal of $T(V)$ contained in $\oplus_{n\geq 2}T^n(V)$. The Nichols algebra of $V$ is
defined by $\Bb(V)=T(V)/I(V)$. $\Bb(V)$ is called diagonal type if 
$(V, c)$ is of digonal type.

By the following lemmas and their proofs, we can find all the irreducible Yetter-Drinfeld modules $M(g, V)$
 in $_G^G\YD$, once we have known the corresponding irreducible representations $(\rho,V)$ of  $\kk G^{g}$. The
 corresponding Nichols algebra of $M(\Oo_g, V)$ is denoted by $\Bb(\Oo_{g}, \rho)$ or $\Bb(\Oo_{g}, V)$.
 \begin{lem}\cite[Lemma 1.4.16]{H1}\label{lem.1.4.16}
 Let $g\in G$,  $M\in {_G^G}\YD(\Oo_g)$. Then $M(g, M_g)\to M$ is an isomorphism of Yetter-Drinfeld modules in $_G^G\YD$.
 \end{lem}
\begin{lem}\cite[Corollary 1.4.18]{H1}\label{cor.1.4.18}
Let $\{\Oo_{g_l}|l\in L\}$ be the set of the conjugacy classes of $G$. There is a bijection between the disjoint union of the isomorphism 
classes of the simple left $\kk G^{g_l}$-modules, $l\in L$, and the set of
isomorphism classes of the simple Yetter-Drinfeld modules in $_G^G\YD$.
\end{lem}

\subsection {The quaternion group $\mathbb{Q}_8$}
As we are familiar,
 $\mathbb{Q}_8=\langle  x,y|x^4=1,y^2=x^2, xy=yx^{-1}\rangle $ is generated by two elements, precisely,
\begin{align*}
\mathbb{Q}_8=\{1,x,x^2,x^3,y,xy,x^2y,x^3y\}
\end{align*}
where $x^{-1}=x^3$, $x^{-2}=x^2$,$y^{-1}=x^2y$, and $(xy)^{-1}=x^3y$.

There are five conjugacy classes in $G=\mathbb{Q}_8$:
\begin{align*}
\Oo_1=\{1\},
\Oo_x=\{x,x^3\},
\Oo_{x^2}=\{x^2\},
\Oo_{y}=\{y,x^2y\},
\Oo_{xy}=\{xy, x^3y\}.
\end{align*}
Choose one element in each conjugacy class, and we compute the centralizers:
\begin{align*}
G^1=G,\quad
G^x=\langle x\rangle\cong \ZZ_4,\quad
G^{x^2}=G,\quad
G^{y}=\langle y\rangle\cong \ZZ_4,\quad
G^{xy}=\langle xy\rangle\cong \ZZ_4.
\end{align*}
By \cite{J1},there are four 1-dimensional irreducible representations $(\rho_i,V_i)_{1\leq i\leq4}\in \mathrm{Irr}(\mathbb{Q}_8)$, and one 2-dimensional irreducible representation, denoted by $(\rho_5,V_5)$. The character table of $(\rho_i,V_i)$ is as follows: \\
\newline
\begin{tabular}{rrrrrrrrrrr}
\hline
    &   1& $x^2$&$x$&$y$&$xy$\\
\hline
$\chi_1$       &   1&1&1&1&1\\
$\chi_2$& 1&1&1&-1&-1\\
$\chi_3$& 1&1&-1&1&-1\\
$\chi_4$& 1&1&-1&-1&1\\
$\chi_5$& 2&-2&0&0&0\\
\hline
\end{tabular}
\subsection{Gelfand-Kirillov dimension}
The Gelfand-Kirillov dimension, GK-dimension for short, 
becomes a powerful tool to study  noncommutative algebras,
especially for  those with infinite dimensions. 
For the definition and properties of the GK-dimension we refer to \cite{G1}.
\begin{defn}\cite{G1}
The Gelfand-Kirillov dimension of a $\kk$-algebra $A$ is 
\begin{align*}
\Gkdim(A)=\sup_V\varlimsup\log_n \dim(V^n),
\end{align*}
where the supremum is taken over all finite dimensional subspaces $V$
of $A$.
\end{defn}
As for finite GK-dimension, we need the following result
\begin{lem}\cite[Theorem 6]{A1}\label{lem2.6}
If either its Weyl groupoid is infinite and $\dim V=2$, or
else $V$ is of affine Cartan type, then $\Gkdim \Bb(V)=\infty$.
\end{lem}

\section{\bf The Nichols algebras $\Bb(\Oo_1,\rho)$}
Consider the conjugacy $\Oo_1=\{1\}$.
$G^1=G=\mathbb{Q}_8$.
For $(\rho_1,V_1)$, the module structure is 
\begin{align*}
a\cdot x_1=x_1, \quad b\cdot x_1=x_1.
\end{align*}
The corresponding Yetter-Drinfeld module $W_1=1\otimes V_1$
The module structure is 
\begin{align*}
a\cdot (1\otimes x_1)=1\otimes \rho_1(a)(x_1)=1\otimes x_1
\end{align*}
The comodule structure is 
\begin{align*}
\rho(1\otimes x_1)=1\otimes (1\otimes x_1)
\end{align*}
Let $w_1=1\otimes x_1$
The braiding structure is 
\begin{align*}
c(w_1\otimes w_1)=w_1\otimes w_1
\end{align*}
Therefore $\dim \Bb(W_1)=\infty$, and $\Gkdim \Bb(W_1)=1$.

For $(\rho_2,V_2)$, the module structure is 
\begin{align*}
a\cdot x_2=x_2, \quad b\cdot x_2=-x_2.
\end{align*}
The corresponding Yetter-Drinfeld module is $W_2=1\otimes V_2$.
The $\kk G$-module structure of $W_2$ is 
\begin{align*}
a\cdot (1\otimes x_2)&=1\otimes \rho_2(a)(x_2)=1\otimes x_2\\
b\cdot (1\otimes x_2)&=1\otimes \rho_2(b)(x_2)=-1\otimes x_2
\end{align*}
The comodule structure is 
\begin{align*}
\rho(1\otimes x_2)=1\otimes (1\otimes x_2)
\end{align*}
Let $w_2=1\otimes x_2$
The braiding structure is 
\begin{align*}
c(w_2\otimes w_2)=w_2\otimes w_2
\end{align*}
Therefore $\dim \Bb(W_2)=\infty$, and $\Gkdim \Bb(W_2)=1$.

For $(\rho_3,V_3)$, the module structure is 
\begin{align*}
a\cdot x_3=-x_3, \quad b\cdot x_3=x_3.
\end{align*}
The corresponding Yetter-Drinfeld module $W_3=1\otimes V_3$
The module structure is 
\begin{align*}
a\cdot (1\otimes x_3)&=1\otimes \rho_3(a)(x_3)=-1\otimes x_3\\
b\cdot (1\otimes x_3)&=1\otimes \rho_3(b)(x_3)=1\otimes x_3
\end{align*}
The comodule structure is 
\begin{align*}
\rho(1\otimes x_3)=1\otimes (1\otimes x_3)
\end{align*}
Let $w_3=1\otimes x_3$
The braiding structure is 
\begin{align*}
c(w_3\otimes w_3)=w_3\otimes w_3
\end{align*}
Therefore $\dim \Bb(W_3)=\infty$, and $\Gkdim \Bb(W_3)=1$.

For $(\rho_4,V_4)$, the module structure is 
\begin{align*}
a\cdot x_4=-x_4, \quad b\cdot x_4=-x_4.
\end{align*}
The corresponding Yetter-Drinfeld module $W_4=1\otimes V_4$
The module structure is 
\begin{align*}
a\cdot (1\otimes x_4)=1\otimes \rho_4(a)(x_4)=-1\otimes x_4\\
b\cdot (1\otimes x_4)=1\otimes \rho_4(b)(x_3)=-1\otimes x_4
\end{align*}
The comodule structure is 
\begin{align*}
\rho(1\otimes x_4)=1\otimes (1\otimes x_4)
\end{align*}
Let $w_4=1\otimes x_4$
The braiding structure is 
\begin{align*}
c(w_4\otimes w_4)=w_4\otimes w_4
\end{align*}
Therefore $\dim \Bb(W_4)=\infty$, and $\Gkdim \Bb(W_4)=1$.

For $(\rho_5,V_5)$, let $\{v_1,v_2\}$ be a basis of $V_5$. The module structure is 
\begin{align*}
a\cdot v_1&=\mathbf{i}v_1,\quad  b\cdot v_1=v_2,\quad
a\cdot v_2=-\mathbf{i}v_2,\quad   b\cdot v_2=-v_1.
\end{align*}
The corresponding Yetter-Drinfeld module $W_5=1\otimes V_5$
The module structure is 
\begin{align*}
a\cdot (1\otimes v_1)&=1\otimes \rho_5(a)(v_1)=\mathbf{i}\otimes v_1,\\
b\cdot (1\otimes v_1)&=1\otimes \rho_5(b)(v_1)=1\otimes v_2,\\
a\cdot (1\otimes v_2)&=1\otimes \rho_5(a)(v_2)=-\mathbf{i}(1\otimes v_2),\\
b\cdot (1\otimes v_2)&=1\otimes \rho_5(b)(v_2)=-1\otimes v_1.
\end{align*}
The comodule structure is 
\begin{align*}
\rho(1\otimes v_1)=1\otimes (1\otimes v_1),\quad
\rho(1\otimes v_2)=1\otimes (1\otimes v_2).
\end{align*}
Let $u_1=1\otimes v_1$, $u_2=1\otimes v_2$.
The braiding structure is 
\begin{align*}
c(u_1\otimes u_1)&=u_1\otimes u_1,\quad
c(u_1\otimes u_2)=u_2\otimes u_1,\\
c(u_2\otimes u_1)&=u_1\otimes u_2,\quad
c(u_2\otimes u_2)=u_2\otimes u_2.
\end{align*}
Therefore $\dim \Bb(W_4)=\infty$, and $\Gkdim \Bb(W_4)=2$ by \cite[Example 31]{A1}.

Therefore, we have proved
\begin{prop}\label{prop:3.1}
For any irreducible representation 
$(\rho, V)\in \Irr (G^{1})$,  $\Gkdim \Bb(\Oo_{1}, \rho)<\infty$, and $\dim \Bb(\Oo_{1}, \rho)=\infty$.
\end{prop}

\section{\bf The Nichols algebras $\Bb(\Oo_{x},\phi)$}
As shown above, $\Oo_{x}=\{x,x^3\}$.
$G^x=\{1,x,x^2,x^3\}\cong \ZZ_4$ is the cyclic group with order 4. Therefore,  all the irrreducible representations are $(\phi_t,U_t), 0\leq t\leq 3$, which are all 1-dimensional. Precisely, 
\begin{align*}
x^k\cdot u_t=\phi_t(x^k)(u_t)=e^{t\pi \mathbf{i}k/2}u_t.
\end{align*}
where $\mathbf{i}=\sqrt{-1}$.

Clearly, $(\phi_0,U_0)$ is the trivial representation of $G^x$.
Since 
\begin{align*}
G=G^x\cup yG^x,
\end{align*}
the corresponding Yetter-Drinfeld module is 
$M(\Oo_{x},\phi_0)=\kk  G\otimes_{\mathbbm{k}G^x}U_0 =1\otimes U_0\oplus y\otimes U_0$.

Let $v_1=1\otimes u_0$, $v_2=y\otimes u_0$.
The $G$-module structure is 
\begin{align*}
x\cdot v_1&=x\cdot (1\otimes u_0)=1\otimes \phi_0(x)(u_0)=1\otimes u_0=v_1,\\
y\cdot v_1&=y\cdot (1\otimes u_0)=y\otimes \phi_0(1)(u_0)=y\otimes u_0=v_2.\\
x\cdot v_2&=x\cdot (y\otimes u_0)=y\otimes \phi_0(x^3)(u_0)
=y\otimes u_0=v_2,\\
y\cdot v_2&=y\cdot (y\otimes u_0)=1\otimes \phi_0(x^2)(u_0)
=1\otimes u_0=v_1.
\end{align*}

The $G$-comdule is 
\begin{align*}
\delta(v_1)&=\delta(1\otimes u_0)=(1\rhd x)\otimes v_1=x\otimes v_1,\\
\delta(v_2)&=\delta(y\otimes u_0)=(y\rhd x)\otimes v_2=x^3\otimes v_2.
\end{align*}
The braiding structure is 
\begin{align*}
c(v_1\otimes v_1)&=v_1^{(-1)}\cdot v_1\otimes v_1^{(0)}
=x\cdot v_1\otimes v_1=v_1\otimes v_1,\\
c(v_1\otimes v_2)&=v_1^{(-1)}\cdot v_2\otimes v_1^{(0)}
=x\cdot v_2\otimes v_1=v_2\otimes v_1,\\
c(v_2\otimes v_1)&=v_2^{(-1)}\cdot v_1\otimes v_2^{(0)}
=x^3\cdot v_1\otimes v_2=v_1\otimes v_2,\\
c(v_2\otimes v_2)&=v_2^{(-1)}\cdot v_2\otimes v_2^{(0)}
=x^3\cdot v_2\otimes v_2=v_2\otimes v_2,
\end{align*}
which is of diagonal type with braiding matrix 
\begin{align*}
\mathbf{q}=
\left[
\begin{array}{ccccccc}
1&1\\
1&1
\end{array}
\right]
\end{align*}
Therefore, $\dim \Bb(\Oo_{x},\phi_0)=\infty$, and $\Gkdim \Bb(\Oo_{x},\phi_0)=2$ by \cite[Example 31]{A1}.

For $(\phi_1,U_1)$, we compute the Yetter-Drinfeld module
\begin{align*}
M(\Oo_x, \phi_1)=1\otimes U_1\oplus y\otimes U_1.
\end{align*}
Let $v_1=1\otimes u_1$, $v_2=y\otimes u_1$. $q=e^{\pi\mathbf{i}/2}\in \GG_{4}'$.The $G$-module
structure is 
\begin{align*}
x\cdot v_1&=x\cdot (1\otimes u_1)=1\otimes \phi_1(x)(u_1)=e^{\pi\mathbf{i}/2}v_1,\\
y\cdot v_1&=y\cdot (1\otimes u_1)=y\otimes \phi_1(1)(u_1)
=v_2.
\end{align*}
The $G$-comodule is 
\begin{align*}
\delta(v_1)&=\delta(1\otimes u_1)=(1\rhd x)\otimes v_1=x\otimes v_1,\\
\delta(v_2)&=\delta(y\otimes u_1)=(y\rhd x)\otimes v_2=x^3\otimes v_2.
\end{align*}
The braiding structure is 
\begin{align*}
c(v_1\otimes v_1)&=v_1^{(-1)}\cdot v_1\otimes v_1^{(0)}
=x\cdot v_1\otimes v_1=qv_1\otimes v_1,\\
c(v_1\otimes v_2)&=v_1^{(-1)}\cdot v_2\otimes v_1^{(0)}
=x\cdot v_2\otimes v_1=q^3v_2\otimes v_1,\\
c(v_2\otimes v_1)&=v_2^{(-1)}\cdot v_1\otimes v_2^{(0)}
=x^3\cdot v_1\otimes v_2=q^3v_1\otimes v_2,\\
c(v_2\otimes v_2)&=v_2^{(-1)}\cdot v_2\otimes v_2^{(0)}
=x^3\cdot v_2\otimes v_2=qv_2\otimes v_2.
\end{align*}

This is of diagonal type with braiding matrix
\begin{align*}
\mathbf{q}=
\left[
\begin{array}{ccccccc}
q&q^3\\
q^3&q
\end{array}
\right]
\end{align*}
which is of affine Cartan type. Therefore $\dim \Bb(\Oo_x, \phi_1)=\infty$, and  $\Gkdim \Bb(\Oo_x, \phi_1)=\infty$ by going through the
list of Heckenberger's classification \cite{H1} and Lemma \ref{lem2.6}.

For the $(\phi_2,U_2)$,
 we compute the Yetter-Drinfeld module
\begin{align*}
M(\Oo_x, \phi_2)=1\otimes U_2\oplus y\otimes U_2.
\end{align*}
Let $v_1=1\otimes u_2$, $v_2=y\otimes u_2$. $q=e^{\pi\mathbf{i}/2}\in \GG_{4}'$.The $G$-module
structure is 
\begin{align*}
x\cdot v_1&=x\cdot (1\otimes u_2)=1\otimes \phi_2(x)(u_2)=e^{2\pi\mathbf{i}/2}(1\otimes u_2)=q^2v_1,\\
y\cdot v_1&=y\cdot (1\otimes u_2)=y\otimes \phi_2(1)(u_2)=v_2,\\
x\cdot v_2&=x\cdot (y\otimes u_2)=y\otimes \phi_2(x^3)(u_2)=q^2v_2,\\
y\cdot v_2&=y\cdot (y\otimes u_2)=1\otimes \phi_2(x^2)(u_2)=q^2v_1.
\end{align*}
The $\kk G$-comodule is 
\begin{align*}
\delta(v_1)&=\delta(1\otimes u_2)=(1\rhd x)\otimes v_1=x\otimes v_1,\\
\delta(v_2)&=\delta(y\otimes u_2)=(y\rhd x)\otimes v_2=x^3\otimes v_2.
\end{align*}
The braiding structure is 
\begin{align*}
c(v_1\otimes v_1)&=v_1^{(-1)}\cdot v_1\otimes v_1^{(0)}
=x\cdot v_1\otimes v_1=q^2v_1\otimes v_1,\\
c(v_1\otimes v_2)&=v_1^{(-1)}\cdot v_2\otimes v_1^{(0)}
=x\cdot v_2\otimes v_1=q^2v_2\otimes v_1,\\
c(v_2\otimes v_1)&=v_2^{(-1)}\cdot v_1\otimes v_2^{(0)}
=x^3\cdot v_1\otimes v_2=q^2v_1\otimes v_2,\\
c(v_2\otimes v_2)&=v_2^{(-1)}\cdot v_2\otimes v_2^{(0)}
=x^3\cdot v_2\otimes v_2=q^2v_2\otimes v_2,\\
\end{align*}

This is of diagonal type with braiding matrix
\begin{align*}
\mathbf{q}=
\left[
\begin{array}{ccccccc}
-1&-1\\
-1&-1
\end{array}
\right]
\end{align*}
By \cite[Example 27]{A1}, $\dim \Bb(\Oo_x, \phi_2)=4$, and  $\Gkdim \Bb(\Oo_x, \phi_2)=0$.


For the $(\phi_3,U_3)$,
 we compute the Yetter-Drinfeld module
\begin{align*}
M(\Oo_x, \phi_3)=1\otimes U_3\oplus y\otimes U_3.
\end{align*}
Let $v_1=1\otimes u_3$, $v_2=y\otimes u_3$. $q=e^{\pi\mathbf{i}/2}\in \GG_{4}'$. The $G$-module
structure is 
\begin{align*}
x\cdot v_1&=x\cdot (1\otimes u_3)=1\otimes \phi_1(x)(u_3)=e^{3\pi\mathbf{i}/2}(1\otimes u_3)=e^{3\pi\mathbf{i}/2}v_1,\\
y\cdot v_1&=y\cdot (1\otimes u_3)=y\otimes \phi_1(1)(u_3)
=v_2.
\end{align*}
The $G$-comodule is 
\begin{align*}
\delta(v_1)&=\delta(1\otimes u_3)=(1\rhd x)\otimes v_1=x\otimes v_1,\\
\delta(v_2)&=\delta(y\otimes u_3)=(y\rhd x)\otimes v_2=x^3\otimes v_2.
\end{align*}
The braiding structure is 
\begin{align*}
c(v_1\otimes v_1)&=v_1^{(-1)}\cdot v_1\otimes v_1^{(0)}
=x\cdot v_1\otimes v_1=e^{3\pi\mathbf{i}/2}v_1\otimes v_1,\\
c(v_1\otimes v_2)&=v_1^{(-1)}\cdot v_2\otimes v_1^{(0)}
=x\cdot v_2\otimes v_1=e^{3\pi\mathbf{i}/2}v_2\otimes v_1,\\
c(v_2\otimes v_1)&=v_2^{(-1)}\cdot v_1\otimes v_2^{(0)}
=x^3\cdot v_1\otimes v_2=e^{\pi\mathbf{i}/2}v_1\otimes v_2,\\
c(v_2\otimes v_2)&=v_2^{(-1)}\cdot v_2\otimes v_2^{(0)}
=x^3\cdot v_2\otimes v_2=e^{3\pi\mathbf{i}/2}v_2\otimes v_2,\\
\end{align*}

This is of diagonal type with braiding matrix
\begin{align*}
\mathbf{q}=
\left[
\begin{array}{ccccccc}
e^{3\pi\mathbf{i}/2}&e^{3\pi\mathbf{i}/2}\\
e^{\pi\mathbf{i}/2}&e^{\pi\mathbf{i}/2}
\end{array}
\right]
\end{align*}
By \cite[Example 27]{A1}, $\dim \Bb(\Oo_x, \phi_3)=4$, and  $\Gkdim \Bb(\Oo_x, \phi_3)=0$.

\begin{prop}\label{prop:4.1}
For any  irreducible representation 
$(\phi, V)\in \Irr (G^{x})$,  $\Gkdim \Bb(\Oo_{x}, \phi)<\infty$ if and only if $\phi=\phi_0$  or $\phi=\phi_2$;
$\dim \Bb(\Oo_{x}, \phi)<\infty$ if and only if
 $\phi=\phi_2$.
\end{prop}

\section{\bf The Nichols algebras $\Bb(\Oo_{x^2},\rho)$}
It is shown that $\Oo_{x^2}=\{x^2\}$, and $G^{x^2}=G$.
It is shown that 
there are four 1-dim irreducible representations $(\rho_1,V_1),(\rho_2,V_2),(\rho_3,V_3),(\rho_4,V_4)$ and one 2-dim 
irreducible representation, denoted by $(\rho_5,V_5)$.

For $(\rho_1,V_1)$
The corresponding Yetter-Drinfeld module $M(\Oo_{x^2},\rho_1)=1\otimes V_1$.
Let $w_1=1\otimes v_1$.
The module structure is 
\begin{align*}
x\cdot w_1&=x\cdot (1\otimes v_1)=1\otimes \rho_1(x)(v_1)=1\otimes v_1=w_1,\\
y\cdot w_1&=y\cdot (1\otimes v_1)=1\otimes \rho_1(y)(v_1)=1\otimes v_1=w_1.
\end{align*}
The comodule structure is 
\begin{align*}
\rho(w_1)=\rho(1\otimes v_1)=x^2\otimes (1\otimes v_1)=x^2\otimes w_1.
\end{align*}

The braiding structure is
\begin{align*}
c(w_1\otimes w_1)=w_1^{(-1)}\cdot w_1\otimes w_1^{(0)}
=x^2\cdot w_1\otimes w_1=w_1\otimes w_1.
\end{align*}

Therefore, $\dim \Bb(\Oo_{x^2},\rho_1)=\infty$, and $\Gkdim \Bb(\Oo_{x^2},\rho_1)=1$.

For $(\rho_2,V_2)$
The corresponding Yetter-Drinfeld module $M(\Oo_{x^2},\rho_2)=1\otimes V_2$.
Let $w_2=1\otimes v_2$.
The module structure is 
\begin{align*}
x\cdot w_2&=x\cdot (1\otimes v_2)=1\otimes \rho_2(x)(v_2)=1\otimes v_2=w_2,\\
y\cdot w_2&=y\cdot (1\otimes v_2)=1\otimes \rho_2(y)(v_2)=-1\otimes v_2=-w_2.
\end{align*}
The comodule structure is 
\begin{align*}
\rho(1\otimes v_2)=x^2\otimes (1\otimes v_2)=x^2\otimes w_2.
\end{align*}

The braiding structure is
\begin{align*}
c(w_2\otimes w_2)=w_2^{(-1)}\cdot w_2\otimes w_2^{(0)}
=x^2\cdot w_2\otimes w_2=w_2\otimes w_2.
\end{align*}

Therefore, $\dim \Bb(\Oo_{x^2},\rho_2)=\infty$, and $\Gkdim \Bb(\Oo_{x^2},\rho_2)=1$.


For $(\rho_3,V_3)$
The corresponding Yetter-Drinfeld module $M(\Oo_{x^2},\rho_3)=1\otimes V_3$.
Let $w_3=1\otimes v_3$.
The module structure is 
\begin{align*}
x\cdot w_3&=x\cdot (1\otimes v_3)=1\otimes \rho_3(x)(v_3)=-1\otimes v_3=-w_3,\\
y\cdot w_3&=y\cdot (1\otimes v_3)=1\otimes \rho_3(y)(v_3)=1\otimes v_3=w_3.
\end{align*}
The comodule structure is 
\begin{align*}
\rho(1\otimes v_3)&=x^2\otimes (1\otimes v_3)=x^2\otimes w_3.
\end{align*}

The braiding structure is
\begin{align*}
c(w_3\otimes w_3)=w_3^{(-1)}\cdot w_3\otimes w_3^{(0)}
=x^2\cdot w_3\otimes w_3=w_3\otimes w_3.
\end{align*}

Therefore, $\dim \Bb(\Oo_{x^2},\rho_3)=\infty$, and $\Gkdim \Bb(\Oo_{x^2},\rho_3)=1$.

For $(\rho_4,V_4)$
The corresponding Yetter-Drinfeld module $M(\Oo_{x^2},\rho_4)=1\otimes V_4$.
Let $w_4=1\otimes v_4$.
The module structure is 
\begin{align*}
x\cdot w_4&=x\cdot (1\otimes v_4)=1\otimes \rho_4(x)(v_4)=-1\otimes v_4=-w_4,\\
y\cdot w_4&=y\cdot (1\otimes v_4)=1\otimes \rho_4(y)(v_4)=1\otimes v_4=w_4.
\end{align*}
The comodule structure is 
\begin{align*}
\rho(1\otimes v_4)=x^2\otimes (1\otimes v_4)=x^2\otimes w_4.
\end{align*}

The braiding structure is
\begin{align*}
c(w_4\otimes w_4)=w_4^{(-1)}\cdot w_4\otimes w_4^{(0)}
=x^2\cdot w_4\otimes w_4=w_4\otimes w_4.
\end{align*}

Therefore, $\dim \Bb(\Oo_{x^2},\rho_4)=\infty$, and $\Gkdim \Bb(\Oo_{x^2},\rho_4)=1$.

For $(\rho_5,V_5)$, the corresponding Yetter-Drinfeld module 
$M(\Oo_{x^2},\rho_5)=1\otimes V_5$. Let $w_5=1\otimes v_5$, $w_6=1\otimes v_6$. The $G$-module structure is
\begin{align*}
x\cdot w_5&=x\cdot (1\otimes v_5)= 1\otimes \rho_5(x)(v_5)
=\mathbf{i}(1\otimes v_5)=\mathbf{i}w_5,\\
y\cdot w_5&=y\cdot (1\otimes v_5)=1\otimes \rho_5(y)(v_5)
=1\otimes v_6=w_6,\\
x\cdot w_6&=x\cdot (1\otimes v_6)=1\otimes \rho_5(x)(v_6)
=-\mathbf{i}(1\otimes v_6)=-\mathbf{i}w_6,\\
y\cdot w_6&=y\cdot (1\otimes v_6)=1\otimes \rho_5(y)(v_6)
=-1\otimes v_5=-w_5. 
\end{align*}

The $G$-comodule structure is
\begin{align*}
\delta(w_5)=x^2\otimes w_5,\quad
\delta(w_6)=x^2\otimes w_6.
\end{align*}
The braiding structure is 
\begin{align*}
c(w_5\otimes w_5)=w_5^{(-1)}\cdot w_5\otimes w_5^{(0)}
=x^2\cdot w_5\otimes w_5=-w_5\otimes w_5,\\
c(w_5\otimes w_6)=w_5^{(-1)}\cdot w_6\otimes w_5^{(0)}
=x^2\cdot w_6\otimes w_5=-w_6\otimes w_5,\\
c(w_6\otimes w_5)=w_6^{(-1)}\cdot w_5\otimes w_6^{(0)}
=x^2\cdot w_6\otimes w_5=-w_6\otimes w_6,\\
c(w_6\otimes w_6)=w_6^{(-1)}\cdot w_6\otimes w_6^{(0)}
=x^2\cdot w_6\otimes w_6=-w_6\otimes w_6.
\end{align*}
which is of diagonal type with braiding matrix 
\begin{align*}
\mathbf{q}=
\left[
\begin{array}{ccccccc}
-1&-1\\
-1&-1
\end{array}
\right]
\end{align*}
Therefore, $\Bb(\Oo_{x^2},\rho_5)\cong \Lambda M(\Oo_{x^2},\rho_5)$, the exterior algebra over $M(\Oo_{x^2},\rho_5)$, which implies $\dim \Bb(\Oo_{x^2},\rho_5)=4$ and $\Gkdim \Bb(\Oo_{x^2},\rho_5)=0$.

\begin{prop}\label{prop:5.1}
Let $(\rho, V)\in \Irr (G^{x^2})$ be any irreducible representation.\\
 Then 
 $\Gkdim \Bb(\Oo_{x^2}, \rho)<\infty$.
Moreover,  $\dim \Bb(\Oo_{x^2}, \rho)<\infty$ if and only if $\rho=\rho_5$.
\end{prop}


\section{\bf The Nichols algebras $\Bb(\Oo_{y},\phi)$}

For the conjugacy class $\Oo_{y}$, $G^{y}=\{1,y,x^2y, x^2\}=\langle  y\rangle\cong \ZZ_4$. As we known, 
Therefore,  all the irrreducible representations are $(\phi_t,U_t), 0\leq t\leq 3$, which are all 1-dimensional. Precisely, 
\begin{align*}
y^k\cdot u_t&=\phi_t(y^k)(u_t)=e^{t\pi \mathbf{i}k/2}u_t.
\end{align*}
where $\mathbf{i}=\sqrt{-1}$.

Clearly, $(\phi_0,U_0)$ is the trivial representation of $G^{y}$.
Since 
\begin{align*}
G=G^y\cup xG^y,
\end{align*}
the corresponding Yetter-Drinfeld module is 
$M(\Oo_{y},\phi_0)=\kk  G\otimes_{\mathbbm{k}G^y}U_0 =1\otimes U_0\oplus x\otimes U_0$.

Let $v_1=1\otimes u_0$, $v_2=x\otimes u_0$.
The $G$-module structure is 
\begin{align*}
x\cdot v_1&=x\cdot (1\otimes u_0)=x\otimes \phi_0(1)(u_0)=x\otimes u_0=v_2,\\
y\cdot v_1&=y\cdot (1\otimes u_0)=1\otimes \phi_0(y)(u_0)=1\otimes u_0=v_1,\\
x\cdot v_2&=x\cdot (x\otimes u_0)=1\otimes \phi_0(y^2)(u_0)
=1\otimes u_0=v_1,\\
y\cdot v_2&=y\cdot (x\otimes u_0)=x\otimes \phi_0(y^3)(u_0)
=x\otimes u_0=v_2.
\end{align*}

The $G$-comdule is 
\begin{align*}
\delta(v_1)&=\delta(1\otimes u_0)=(1\rhd y)\otimes v_1=y\otimes v_1,\\
\delta(v_2)&=\delta(x\otimes u_0)=(x\rhd y)\otimes v_2=y^3\otimes v_2.
\end{align*}
The braiding structure is 
\begin{align*}
c(v_1\otimes v_1)&=v_1^{(-1)}\cdot v_1\otimes v_1^{(0)}
=y\cdot v_1\otimes v_1=v_1\otimes v_1,\\
c(v_1\otimes v_2)&=v_1^{(-1)}\cdot v_2\otimes v_1^{(0)}
=y\cdot v_2\otimes v_1=v_2\otimes v_1,\\
c(v_2\otimes v_1)&=v_2^{(-1)}\cdot v_1\otimes v_2^{(0)}
=y^3\cdot v_1\otimes v_2=v_1\otimes v_2,\\
c(v_2\otimes v_2)&=v_2^{(-1)}\cdot v_2\otimes v_2^{(0)}
=y^3\cdot v_2\otimes v_2=v_2\otimes v_2,
\end{align*}
which is of diagonal type with braiding matrix 
\begin{align*}
\mathbf{q}=
\left[
\begin{array}{ccccccc}
1&1\\
1&1
\end{array}
\right]
\end{align*}
Therefore, $\dim \Bb(\Oo_{y},\phi_0)=\infty$, and $\Gkdim \Bb(\Oo_{y},\phi_0)=2$ by \cite[Example 31]{A1}.

For  $(\phi_1,U_1)$, 
the corresponding Yetter-Drinfeld module is 
$M(\Oo_{y},\phi_1)=\kk  G\otimes_{\mathbbm{k}G^y}U_1 =1\otimes U_1\oplus x\otimes U_1$.

Let $v_1=1\otimes u_1$, $v_2=x\otimes u_1$.
The $G$-module structure is 
\begin{align*}
x\cdot v_1&=x\cdot (1\otimes u_1)=x\otimes \phi_1(1)(u_1)=x\otimes u_1=v_2,\\
y\cdot v_1&=y\cdot (1\otimes u_1)=1\otimes \phi_1(y)(u_1)=q(1\otimes u_1)=qv_1,\\
x\cdot v_2&=x\cdot (x\otimes u_1)=1\otimes \phi_1(y^2)(u_1)
=q^2(1\otimes u_1)=q^2v_1,\\
y\cdot v_2&=y\cdot (x\otimes u_0)=x\otimes \phi_1(y^3)(u_1)
=q^3(x\otimes u_1)=q^3v_2.
\end{align*}

The $G$-comdule is 
\begin{align*}
\delta(v_1)&=\delta(1\otimes u_1)=(1\rhd y)\otimes v_1=y\otimes v_1,\\
\delta(v_2)&=\delta(x\otimes u_1)=(x\rhd y)\otimes v_2=y^3\otimes v_2.
\end{align*}
The braiding structure is 
\begin{align*}
c(v_1\otimes v_1)&=v_1^{(-1)}\cdot v_1\otimes v_1^{(0)}
=y\cdot v_1\otimes v_1=qv_1\otimes v_1,\\
c(v_1\otimes v_2)&=v_1^{(-1)}\cdot v_2\otimes v_1^{(0)}
=y\cdot v_2\otimes v_1=q^3v_2\otimes v_1,\\
c(v_2\otimes v_1)&=v_2^{(-1)}\cdot v_1\otimes v_2^{(0)}
=y^3\cdot v_1\otimes v_2=q^3v_1\otimes v_2,\\
c(v_2\otimes v_2)&=v_2^{(-1)}\cdot v_2\otimes v_2^{(0)}
=y^3\cdot v_2\otimes v_2=qv_2\otimes v_2,
\end{align*}
which is of diagonal type with braiding matrix 
\begin{align*}
\mathbf{q}=
\left[
\begin{array}{ccccccc}
q&q^3\\
q^3&q
\end{array}
\right]
\end{align*}
Therefore, by going through the
list of Heckenberger's classification \cite{H1} and Lemma \ref{lem2.6}, $\dim \Bb(\Oo_y, \phi_1)=\infty$,
and $\Gkdim \Bb(\Oo_y, \phi_1)=\infty$.

For  $(\phi_2,U_2)$, 
the corresponding Yetter-Drinfeld module is 
\begin{align*}
M(\Oo_{y},\phi_2)=\kk  G\otimes_{\mathbbm{k}G^y}U_2 =1\otimes U_2\oplus x\otimes U_2.
\end{align*}
Let $v_1=1\otimes u_2$, $v_2=x\otimes u_2$.
The $G$-module structure is 
\begin{align*}
x\cdot v_1&=x\cdot (1\otimes u_2)=x\otimes \phi_1(1)(u_2)=x\otimes u_2=v_2,\\
y\cdot v_1&=y\cdot (1\otimes u_2)=1\otimes \phi_1(y)(u_2)=q^2(1\otimes u_2)=qv_1,\\
x\cdot v_2&=x\cdot (x\otimes u_2)=1\otimes \phi_1(y^2)(u_2)
=(1\otimes u_2)=v_1,\\
y\cdot v_2&=y\cdot (x\otimes u_2)=x\otimes \phi_1(y^3)(u_2)
=q^2(x\otimes u_2)=q^2v_2.
\end{align*}

The $G$-comdule is 
\begin{align*}
\delta(v_1)&=\delta(1\otimes u_2)=(1\rhd y)\otimes v_1=y\otimes v_1,\\
\delta(v_2)&=\delta(x\otimes u_2)=(x\rhd y)\otimes v_2=y^3\otimes v_2.
\end{align*}
The braiding structure is 
\begin{align*}
c(v_1\otimes v_1)&=v_1^{(-1)}\cdot v_1\otimes v_1^{(0)}
=y\cdot v_1\otimes v_1=q^2v_1\otimes v_1,\\
c(v_1\otimes v_2)&=v_1^{(-1)}\cdot v_2\otimes v_1^{(0)}
=y\cdot v_2\otimes v_1=q^2v_2\otimes v_1,\\
c(v_2\otimes v_1)&=v_2^{(-1)}\cdot v_1\otimes v_2^{(0)}
=y^3\cdot v_1\otimes v_2=q^2v_1\otimes v_2,\\
c(v_2\otimes v_2)&=v_2^{(-1)}\cdot v_2\otimes v_2^{(0)}
=y^3\cdot v_2\otimes v_2=q^2v_2\otimes v_2,
\end{align*}
which is of diagonal type with braiding matrix 
\begin{align*}
\mathbf{q}=
\left[
\begin{array}{ccccccc}
-1&-1\\
-1&-1
\end{array}
\right]
\end{align*}
Therefore, by \cite[Example 31]{A1}, $\dim \Bb(\Oo_y, \phi_2)=4$,
and $\Gkdim \Bb(\Oo_y, \phi_2)=0$.

For  $(\phi_3,U_3)$, 
the corresponding Yetter-Drinfeld module is 
\begin{align*}
M(\Oo_{y},\phi_3)=\kk  G\otimes_{\mathbbm{k}G^y}U_3 =1\otimes U_3\oplus x\otimes U_3.
\end{align*}
Let $v_1=1\otimes u_3$, $v_2=x\otimes u_3$.
The $G$-module structure is 
\begin{align*}
x\cdot v_1&=x\cdot (1\otimes u_3)=x\otimes \phi_3(1)(u_3)=x\otimes u_3=v_2,\\
y\cdot v_1&=y\cdot (1\otimes u_2)=1\otimes \phi_3(y)(u_3)=q^3(1\otimes u_3)=q^3v_1,\\
x\cdot v_2&=x\cdot (x\otimes u_3)=1\otimes \phi_3(y^2)(u_3)
=q^2(1\otimes u_3)=q^2v_1,\\
y\cdot v_2&=y\cdot (x\otimes u_3)=x\otimes \phi_3(y^3)(u_3)
=q(x\otimes u_3)=qv_2.
\end{align*}

The $\kk G$-comdule structure is 
\begin{align*}
\delta(v_1)&=\delta(1\otimes u_3)=(1\rhd y)\otimes v_1=y\otimes v_1,\\
\delta(v_2)&=\delta(x\otimes u_3)=(x\rhd y)\otimes v_2=y^3\otimes v_2.
\end{align*}
The braiding structure is 
\begin{align*}
c(v_1\otimes v_1)&=v_1^{(-1)}\cdot v_1\otimes v_1^{(0)}
=y\cdot v_1\otimes v_1=q^3 v_1\otimes v_1,\\
c(v_1\otimes v_2)&=v_1^{(-1)}\cdot v_2\otimes v_1^{(0)}
=y\cdot v_2\otimes v_1=qv_2\otimes v_1,\\
c(v_2\otimes v_1)&=v_2^{(-1)}\cdot v_1\otimes v_2^{(0)}
=y^3\cdot v_1\otimes v_2=qv_1\otimes v_2,\\
c(v_2\otimes v_2)&=v_2^{(-1)}\cdot v_2\otimes v_2^{(0)}
=y^3\cdot v_2\otimes v_2=q^3v_2\otimes v_2,
\end{align*}
which is of diagonal type with braiding matrix 
\begin{align*}
\mathbf{q}=
\left[
\begin{array}{ccccccc}
q^3&q\\
q&q^3
\end{array}
\right]
\end{align*}
Therefore, by \cite[Theorem 6]{A1} and going through the list of Heckenberger's classification, $\dim \Bb(\Oo_y, \phi_3)=\infty$,
and $\Gkdim \Bb(\Oo_y, \phi_3)=\infty$.

\begin{prop}\label{prop:6.1}
For any  irreducible representation 
$(\phi, V)\in \Irr (G^{y})$,  $\Gkdim \Bb(\Oo_{y}, \phi)<\infty$ if and only if $\phi=\phi_0$  or $\phi=\phi_2$;
and $\dim \Bb(\Oo_{y}, \phi)<\infty$ if and only if
 $\phi=\phi_2$;
\end{prop}

\section{\bf The Nichols algebras $\Bb(\Oo_{xy},\phi)$}
For the conjugacy class $\Oo_{xy}=\{xy, x^3y\}$, $G^{xy}=\{1,xy, x^2,x^3y\}=\langle  xy\rangle\cong \ZZ_4$. As we known, 
Therefore,  all the irrreducible representations are $(\phi_t,U_t), 0\leq t\leq 3$, which are all 1-dimensional. Precisely, 
\begin{align*}
(xy)^k\cdot u_t&=\phi_t((xy)^k)(u_t)=e^{t\pi \mathbf{i}k/2}u_t.
\end{align*}
where $\mathbf{i}=\sqrt{-1}$.

Clearly, $(\phi_0,U_0)$ is the trivial representation of $G^{xy}$.
Since 
\begin{align*}
G=G^{xy}\cup yG^{xy},
\end{align*}
the corresponding Yetter-Drinfeld module is 
\begin{align*}
M(\Oo_{xy},\phi_0)=\kk  G\otimes_{\mathbbm{k}G^{xy}}U_0 =1\otimes U_0\oplus y\otimes U_0.
\end{align*}
Let $v_1=1\otimes u_0$, $v_2=y\otimes u_0$.
The $G$-module structure is 
\begin{align*}
x\cdot v_1&=x\cdot (1\otimes u_0)=y\otimes \phi_0(xy)(u_0)=y\otimes u_0=v_2,\\
y\cdot v_1&=y\cdot (1\otimes u_0)=y\otimes \phi_0(1)(u_0)=y\otimes u_0=v_2,\\
x\cdot v_2&=x\cdot (y\otimes u_0)=1\otimes \phi_0(xy)(u_0)
=1\otimes u_0=v_1,\\
y\cdot v_2&=y\cdot (y\otimes u_0)=1\otimes \phi_0((xy)^2)(u_0)
=1\otimes u_0=v_1.
\end{align*}
and we have $(xy)\cdot v_1=v_1,
(xy)\cdot v_2=v_2$.

The $\kk G$-comdule is 
\begin{align*}
\delta(v_1)&=\delta(1\otimes u_0)=(1\rhd xy)\otimes v_1=xy\otimes v_1,\\
\delta(v_2)&=\delta(x\otimes u_0)=(x\rhd xy)\otimes v_2=x^3y\otimes v_2=(xy)^3\otimes v_2.
\end{align*}
The braiding structure is 
\begin{align*}
c(v_1\otimes v_1)&=v_1^{(-1)}\cdot v_1\otimes v_1^{(0)}
=xy\cdot v_1\otimes v_1=v_1\otimes v_1,\\
c(v_1\otimes v_2)&=v_1^{(-1)}\cdot v_2\otimes v_1^{(0)}
=xy\cdot v_2\otimes v_1=v_2\otimes v_1,\\
c(v_2\otimes v_1)&=v_2^{(-1)}\cdot v_1\otimes v_2^{(0)}
=(xy)^3\cdot v_1\otimes v_2=v_1\otimes v_2,\\
c(v_2\otimes v_2)&=v_2^{(-1)}\cdot v_2\otimes v_2^{(0)}
=(xy)^3\cdot v_2\otimes v_2=v_2\otimes v_2,
\end{align*}
which is of diagonal type with braiding matrix 
\begin{align*}
\mathbf{q}=
\left[
\begin{array}{ccccccc}
1&1\\
1&1
\end{array}
\right]
\end{align*}
Therefore $\dim \Bb (\Oo_{xy},\phi_0)=\infty$ and $\Gkdim \Bb (\Oo_{xy},\phi_0)=2$ by \cite[Example 31]{A1}.

For $(\phi_1,U_1)$,
the corresponding Yetter-Drinfeld module is 
\begin{align*}
M(\Oo_{xy},\phi_1)=\kk  G\otimes_{\mathbbm{k}G^{xy}}U_1 =1\otimes U_1\oplus x\otimes U_1.
\end{align*}
Let $v_1=1\otimes u_1$, $v_2=x\otimes u_1$.
The $\kk G$-module structure is 
\begin{align*}
x\cdot v_1&=x\cdot (1\otimes u_1)=y\otimes \phi_1(xy)(u_1)=qy\otimes u_1=qv_2,\\
y\cdot v_1&=y\cdot (1\otimes u_1)=y\otimes \phi_1(1)(u_1)=y\otimes u_1=v_2,\\
x\cdot v_2&=x\cdot (y\otimes u_1)=1\otimes \phi_1(xy)(u_1)
=q(1\otimes u_1)=qv_1,\\
y\cdot v_2&=y\cdot (y\otimes u_1)=1\otimes \phi_1((xy)^2)(u_1)
=q^2(1\otimes u_1)=q^2v_1
\end{align*}
and we have $(xy)\cdot v_1=qv_1, (xy)\cdot v_2=q^3v_2$.
The $\kk G$-comdule is 
\begin{align*}
\delta(v_1)&=\delta(1\otimes u_1)=(1\rhd xy)\otimes v_1=xy\otimes v_1,\\
\delta(v_2)&=\delta(y\otimes u_1)=(y\rhd xy)\otimes v_2=x^3y\otimes v_2=(xy)^3\otimes v_2.
\end{align*}
The braiding structure is 
\begin{align*}
c(v_1\otimes v_1)&=v_1^{(-1)}\cdot v_1\otimes v_1^{(0)}
=xy\cdot v_1\otimes v_1=qv_1\otimes v_1,\\
c(v_1\otimes v_2)&=v_1^{(-1)}\cdot v_2\otimes v_1^{(0)}
=xy\cdot v_2\otimes v_1=q^3v_2\otimes v_1,\\
c(v_2\otimes v_1)&=v_2^{(-1)}\cdot v_1\otimes v_2^{(0)}
=(xy)^3\cdot v_1\otimes v_2=q^3v_1\otimes v_2,\\
c(v_2\otimes v_2)&=v_2^{(-1)}\cdot v_2\otimes v_2^{(0)}
=(xy)^3\cdot v_2\otimes v_2=qv_2\otimes v_2,
\end{align*}
Therefore $M(\Oo_{xy},\phi_1)$ is of diagonal type with braiding matrix
\begin{align*}
\left[  
\begin{array}{cccccc}
q&q^3\\
q^3&q
\end{array}
\right]
\end{align*}
Thus, $\dim \Bb(\Oo_{xy},\phi_1)=\infty$ since it is
of Cartan type. Moreover, $\Gkdim\Bb(\Oo_{xy},\phi_1)=\infty$ by comparing the table by Heckenberger.

For $(\phi_2,U_2)$,
the corresponding Yetter-Drinfeld module is 
\begin{align*}
M(\Oo_{xy},\phi_2)=\kk  G\otimes_{\mathbbm{k}G^{xy}}U_2 =1\otimes U_2\oplus x\otimes U_2.
\end{align*}
Let $v_1=1\otimes u_2$, $v_2=x\otimes u_2$.
The $\kk G$-module structure is 
\begin{align*}
x\cdot v_1&=x\cdot (1\otimes u_1)=y\otimes \phi_1(xy)(u_2)=q^2y\otimes u_1=q^2v_2,\\
y\cdot v_1&=y\cdot (1\otimes u_2)=y\otimes \phi_1(1)(u_2)=y\otimes u_2=v_2,\\
x\cdot v_2&=x\cdot (y\otimes u_2)=1\otimes \phi_1(xy)(u_2)
=q^2(1\otimes u_2)=q^2v_1,\\
y\cdot v_2&=y\cdot (y\otimes u_2)=1\otimes \phi_1((xy)^2)(u_2)
=(1\otimes u_2)=v_1
\end{align*}
and we have $(xy)\cdot v_1=q^2v_1, (xy)\cdot v_2=q^2v_2$.
The $\kk G$-comdule is 
\begin{align*}
\delta(v_1)&=\delta(1\otimes u_2)=(1\rhd xy)\otimes v_1=xy\otimes v_1,\\
\delta(v_2)&=\delta(y\otimes u_2)=(y\rhd xy)\otimes v_2=x^3y\otimes v_2=(xy)^3\otimes v_2.
\end{align*}
The braiding structure is 
\begin{align*}
c(v_1\otimes v_1)&=v_1^{(-1)}\cdot v_1\otimes v_1^{(0)}
=xy\cdot v_1\otimes v_1=q^2v_1\otimes v_1,\\
c(v_1\otimes v_2)&=v_1^{(-1)}\cdot v_2\otimes v_1^{(0)}
=xy\cdot v_2\otimes v_1=q^2v_2\otimes v_1,\\
c(v_2\otimes v_1)&=v_2^{(-1)}\cdot v_1\otimes v_2^{(0)}
=(xy)^3\cdot v_1\otimes v_2=q^2v_1\otimes v_2,\\
c(v_2\otimes v_2)&=v_2^{(-1)}\cdot v_2\otimes v_2^{(0)}
=(xy)^3\cdot v_2\otimes v_2=q^2v_2\otimes v_2,
\end{align*}
Therefore $M(\Oo_{xy},\phi_2)$ is of diagonal type with braiding matrix
\begin{align*}
\left[  
\begin{array}{cccccc}
-1&-1\\
-1&-1
\end{array}
\right]
\end{align*}
Thus, $\dim \Bb(\Oo_{xy},\phi_2)=4$ and
 $\Gkdim\Bb(\Oo_{xy},\phi_2)=0$ by \cite[Example 31]{A1}.  

\begin{prop}\label{prop:7.1}
For any  irreducible representation 
$(\phi, V)\in \Irr (G^{xy})$,  $\Gkdim \Bb(\Oo_{xy}, \phi)<\infty$ if and only if $\phi=\phi_0$  or $\phi=\phi_2$;
and $\dim \Bb(\Oo_{xy}, \phi)<\infty$ if and only if
 $\phi=\phi_2$.
\end{prop}

\section{\bf Acknowledgments}
The author is grateful to Prof. G. X. Liu for insightful conversations on the classification of Nichols algebras.

\bigskip
{\bf Disclosure statement:} The authors report there are no competing interests to declare.


\begin{thebibliography}{99}
\bibitem{A1}Andruskiewitsch N.: An Introduction to Nichols Algebras. In Quantization, Geometry and Noncommutative Structures in Mathematics and Physics. A. Cardona, P. Morales, H. Ocampo, S. Paycha, A. Reyes, eds., Springer, 135-195 (2017)
\bibitem{A2}Andruskiewitsch N., Heckenberger I. and Schneider H. J.: On the Classification of  Finite-dimensional 
Pointed Hopf Algebras. Ann. Math. 171, 375-417 (2010)
\bibitem{A3}Andruskiewitsch N. and  Angiono I.: On Finite
Dimensional Nichols Algebras of Diagonal Type. Bull. Math. Sci.  7, 353-573 (2017)
\bibitem{A4}Andruskiewitsch N., Fantino F. Grana M. and 
Vendramin L.: Pointed Hopf algebras over the Sporadic Simple Groups. J. Algebra. 325,  305-320 (2011)
\bibitem{A5}Andruskiewitsch N., Angiono I. and Heckenberger I.: On Finite GK-Dimensional Nichols Algebras over Abelian
Groups. Mem. Amer. Math. Soc. 2711329 (2021) 
\bibitem{A6}Andruskiewitsch N., Angiono I. and Heckenberger I.: On Nichols Algebras of Infinite Rank with Finite Gelfand-Kirillov Dimension. Rend. Lincei Mat. Appl. 31, 81-101 (2020)
\bibitem{A7}Andruskiewitsch N., Heckenberger I. and Schneider H. J.: The Nichols Algebras of a Semisimple
Yetter-Drinfeld Module. arXiv: 0803.2430v1
\bibitem{A8}Andruskiewitsch N., Carnovale G. and García, G. A.: Finite-dimensional Pointed Hopf Algebras Over Finite Simple Groups of Lie Type I. Non-semisimple Classes in PSL$_n$(q). J. Algebra 442, 36–65 (2015)
\bibitem{A9}Andruskiewitsch, N., Carnovale G. and García G. A.: Finite-dimensional Pointed Hopf Algebras over Finite Simple Groups of Lie type II: Unipotent Classes in Symplectic Groups. Commun. Contemp. Math. 18(4), 1550053, 35 pp (2016)
\bibitem{A10} Andruskiewitsch, N., Carnovale, G. and García, G. A.: Finite-dimensional Pointed Hopf Algebras over Finite Simple Groups of Lie type III. Semisimple Classes in PSL$_n(q)$. Rev. Mat. Iberoam. 33(3), 995–1024 (2017)
\bibitem{A11}Andruskiewitsch, N., Carnovale, G. and García, G. A.: Finite-dimensional Pointed Hopf Algebras over Finite Simple Groups of Lie type IV. Unipotent Classes in Chevalley and Steinberg Groups. Algebr. Represent. Theory 23(3), 621–655 (2020)
\bibitem{A12} Andruskiewitsch, N., Carnovale, G. and  García, G. A.: Finite-dimensional Pointed Hopf Algebras over Finite Simple Groups of Lie type V. Mixed Classes in Chevalley and Steinberg Groups. Manuscripta Math. 166, 605–647  (2021)
\bibitem{A13} Andruskiewitsch, N.: On Pointed Hopf Algebras
over Nilpotent Groups. arXiv: 2104. 04789v1
\bibitem{A}Assem I., Simson D. and Skowro\'{n}ski A.: Elements of Representation Theory of Associative Algebras Volume 1: Techniques of
Representation Theory, London Mathematical Society Student Texts 65, Cambridge University Press, 2006
\bibitem{C}Carnovale, G. and Costantini, M.: Finite-dimensional Pointed Hopf Algebras over Finite Simple Groups of Lie type VI. Suzuki and Ree Groups. J. Pure Appl. Algebra 225(4), 106568, 19 pp (2021)

\bibitem{H0}Heckenberger I.:
The Weyl Groupoid of a Nichols Algebra of Diagonal Type. 
Invent. Math. 164(1), 175–188 (2006)
\bibitem{H1}Heckenberger I. and Schneider H. J.: Hopf
algebras and Root Systems, Mathematatical Surveys and Monographs 247, Amer. Math Soc. (2020)
\bibitem{H1}Heckenberger I.:  Classification of arithmetic root systems. Adv. Math. 220, 59-124(2009)
\bibitem{G1}Krause G. and Lenagan T.: Growth of Algebras and Gelfand-Kirillov Dimension. Revised edition. Graduate Studies in Math. 22. Amer. Math Soc. (2000)
\bibitem{J1} James,G. D. and Liebeak, M. W. Representations and characters of groups. 2 Eds. New York: Cambridge University,
2001.
\bibitem{W1} Wang D. G.,  Zhang J. J. and  Zhuang G.:
Connected Hopf Algebras of Gelfand-Kirillov Dimension Four. Transactions of the American Mathematical Society.
Volume 367(8), 5597–5632  (2015)
\bibitem{Wu} Wu J. Y., Liu G. X. and Ding N. Q.: Classification of Affine Prime Regular Hopf Algebras of GK-dimension One. Adv. Math. 296, 1–54 (2016)
\bibitem{Z0} Zhang Y. L.:Finite GK-Dimensional Nichols Algebras Over the Infinite Dihedral Group. Algebr Represent Theor (2023). https://doi.org/10.1007/s10468-023-10213-1
\bibitem{Z} Zhuang G. B.: Properties of pointed and connected Hopf algebras of finite Gelfand–Kirillov 
dimension. J. London Math. Soc.  87(2) , 877–898 (2013)
\end{thebibliography}
\end{document}